\documentclass[11pt]{article}
\usepackage{amsfonts}
\usepackage{amssymb}

\newtheorem{theorem}{Theorem}[section]

\newtheorem{corollary}[theorem]{Corollary}
\newtheorem{definition}{Definition}[section]
\newtheorem{remark}{Remark}[section]

 \usepackage{amssymb}
\usepackage{epsfig}
\usepackage{amsmath}
\usepackage{amsfonts}

\newcommand{\R}{{{\mathbb R}}}

\def\qed{\hbox to 0pt{}\hfill$\rlap{$\sqcap$}\sqcup$\medbreak}

\title{Greatest solutions and differential inequalities: a journey in two directions \footnote{Partially
supported by 
Ministerio de Educaci\'on y Ciencia, Spain,
project MTM2010-15314, and by Xunta de Galicia, project 2011/XA023.}}
\begin{document}
\maketitle

\vspace{1cm}

\begin{center}
{\large Rodrigo L\'opez Pouso \\
Department of Mathematical Analysis\\
Faculty of Mathematics,\\University of Santiago de Compostela, Campus Sur\\
15782 Santiago de
Compostela, Spain.
{\bf e-mail:} rodrigo.lopez@usc.es
}
\end{center}

\begin{abstract}
We present a new elementary proof of the existence of the least and the greatest solutions to initial value problems in the conditions of Peano's existence theorem. Our proof is based on a modification of Perron's method which allows us to obtain quite easily the greatest solution as the solution with biggest possible integral. In doing so, we simplify the usual proofs, technically overloaded with lower (upper) solutions and/or related differential inequalities. Moreover, those differential (and integral) inequalities, which are interesting in their own right, can be quickly proven by means of known techniques once we know that the greatest and the least solutions exist. Summing up, we revert the usual approach from differential inequalities to extreme solutions, getting a somewhat smoother presentation.\end{abstract}


 \section{Introduction}
Let $(t_0,x_0) \in \R^2$ be fixed and consider the initial value problem
\begin{equation}
\label{ivp}
x'(t)=f(t,x(t)), \, \, t \ge t_0, \, \, \, \, x(t_0)=x_0,
\end{equation}  
where  $f(t,x)$ is real--valued and continuous on a neighborhood of $(t_0,x_0)$.

Under these assumptions, Peano's existence theorem guarantees the existence of at least one solution of (\ref{ivp}). Furthermore, even the existence of the greatest and the least solutions of (\ref{ivp}) can be ensured, and such a result was already proven by Peano in his original paper \cite{peano}! 

As readers can check in \cite{agor, codlev, halanay, lale, reid, sansone, wal}, the available proofs that (\ref{ivp}) has the greatest and the least solution rely either on lower and upper solutions or on differential inequalities which are studied first.  In Section 2 we offer a new concise proof of the existence of the greatest and the least solutions of (\ref{ivp}) which avoids differential inequalities and any appeal to lower and upper solutions. In Section 3, with the existence of the greatest and the least solutions at hand, we proceed to establish classical differential inequalities by means of old yet maybe optimal techniques. In doing so, we get the well--known formula for greatest solutions as suprema of lower solutions one of the typical starting points for proving that greatest solutions exist. Section 4 is devoted to strict differential inequalities and , finally, Section 5 discusses why we shall call greatest solution to what most mathematicians call maximal solution.

We have considered solvability of (\ref{ivp}) on the right of $t_0$ for simplicity. The corresponding assumptions and results on the left of $t_0$ can be obtained through change of variables.
 \section{New proof of existence of greatest solutions}
 This section contains a very easy proof of the existence of the greatest and the least solutions of (\ref{ivp}). Interestingly, our proof does not lean on differential inequalities nor on special sequences of approximate solutions. After any one of the usual proofs of Peano's existence theorem, {\it cf.} \cite[Theorem 2.1]{hart}, the first ideas in the proof are just repetitions.
 
 \bigbreak
 
 Here and henceforth, $D$ is a relatively open neighborhood of $(t_0,x_0)$ in the semiplane $t \ge t_0$.
 
 \begin{theorem}
 \label{pro1}
 If $f:D \longrightarrow \R$ is continuous in $D$ then there exists $L>0$ such that (\ref{ivp}) has the greatest and the least solutions defined on the interval $I=[t_0,t_0+L]$, i.e., there exists a couple of solutions (maybe identical) $x_*, \, x^*:I\longrightarrow \R$ such that any other solution $x:I \longrightarrow \R$ satisfies $$x_*(t) \le x(t) \le x^*(t) \quad \mbox{for all $t \in I$.}$$
 \end{theorem}

\noindent
{\bf Proof.} First, we construct a rectangle $R \subset \R^2$ which contains the graphs of all possible solutions of (\ref{ivp}): let us fix $a>0$, $b>0$ and $M>0$ such that 
\begin{equation}
\label{bded}
|f(t,x)| \le M \quad \mbox{for all $(t,x) \in R=[t_0,t_0+a] \times [x_0-b,x_0+b]$,}
\end{equation}
and let us consider an interval $I=[t_0,t_0+L]$ with length
$$L=\min \left\{a, \dfrac{b}{M+1}\right\}.$$

Peano's existence theorem guarantees that (\ref{ivp}) has at least one solution defined on the interval $I$. Moreover, standard arguments show that every solution $x:I \longrightarrow \R$ satisfies $(t,x(t)) \in \mbox{Int}(R)$ for all $t \in I$, and, as a result, $|x'(t)| = |f(t,x(t))| \le M$ for all $t \in I$. The Ascoli--Arzel\`a Theorem then ensures that ${\cal S}$, the set of solutions of (\ref{ivp}) on the interval $I$, is a relatively compact set of ${\cal C}(I)$, equipped with its usual topology of uniform convergence. Further, the continuity of $f$ implies that ${\cal S}$ is closed, i.e. uniform limits of elements of ${\cal S}$ belong to ${\cal S}$ as well, and therefore ${\cal S}$ is compact. Hence, the continuous mapping
$$x \in {\cal S} \longmapsto \int_{t_0}^{t_0+L}{x(s) \, ds} \in \R,$$
attains a maximum in ${\cal S}$, i.e., there exists at least one $x^* \in {\cal S}$ such that
\begin{equation}
\label{max}
\int_{t_0}^{t_0+L}{x^*(s) \, ds} \ge \int_{t_0}^{t_0+L}{x(s) \, ds} \quad \mbox{for all $x \in {\cal S}$.}
\end{equation}

Let us show that $x^*$ is the greatest solution of (\ref{ivp}) on the interval $I$. Reasoning by contradiction, assume that we have a solution $x:I \longrightarrow \R$ such that $x(t_1)>x^*(t_1)$ for some $t_1 \in (t_0,t_0+L)$. Then the pointwise maximum of $x$ and $x^*$ belongs to ${\cal S}$ and does not satisfy (\ref{max}), a contradiction. 

We can prove in a similar way that the solution having minimum integral is the least solution of (\ref{ivp}). \qed

Knowing that the least and the greatest solutions exist is interesting in its own right and it is also useful as a theoretical tool. For instance, the proof of the following uniqueness result (also due to Peano) cannot be easier.

\begin{corollary}
\label{peaun}
In the conditions of Theorem \ref{pro1}, if $f(t,x)$ is nonincreasing with respect to its second variable, then (\ref{ivp}) has a unique solution.
\end{corollary}

\noindent
{\bf Proof.} Let $L>0$, $x_*$ and $x^*$ be as in Theorem \ref{pro1}. It suffices to prove that $x_*=x^*$ on $I$. To do so, we note that for any $t \in I$ we have
$$x_*(t) \le x^*(t)=x_0+\int_{t_0}^{t}{f(s,x^*(s)) \, ds} \le x_0+\int_{t_0}^{t}{f(s,x_*(s)) \, ds}=x_*(t).$$
\qed
 
In the next section we shall need global information in terms of noncontinuable solutions. It appears, however, that the existence of the greatest (or the least) noncontinuable solution wants a little bit more analysis than expected at first sight.  

\bigbreak

The following variation of Peano's example will highlight the type of difficulties we have when passing from local to global solutions. The initial value problem
\begin{equation}
\label{expepo}
x'=\left\{
\begin{array}{cl}
3x^{2/3}, & \mbox{if $x \le 1$,} \\
\\
3x^2, & \mbox{if $x >1$,}
\end{array}
\right\}, \, \, t \ge 0, \quad x(0)=0,
\end{equation}
has, for each $\tau \ge 0$, a solution given by 
$$\varphi_{\tau}(t)=\left\{
\begin{array}{cl}
0, & \mbox{if $t�\in [0,\tau]$,} \\
(t-\tau)^3, & \mbox{if $t \in [\tau, \tau+1]$,}\\
\dfrac{1}{1-3(t-\tau-1)}, & \mbox{if $t \in [\tau+1,\tau+4/3)$.}
\end{array}
\right.
$$
This justifies the following claims:
\begin{enumerate}
\item We may have infinitely many solutions on a given interval $[0,T]$ and not have the greatest solution {\it on that interval}.
\item We may have different greatest solutions on different intervals: namely, $\varphi_\tau$ is the greatest solution on the interval $[0,\tau+4/3)$.
\item The ``natural" greatest solution of (\ref{expepo}) is $\varphi_0$, and the least solution is the zero function. Notice that their respective intervals of definition are different.
\end{enumerate}

\bigbreak

Recall that if $x:I_1 \longrightarrow \R$ is a solution of (\ref{ivp}), then a continuation of $x$ is another solution $y:I_2\longrightarrow \R$ such that $I_1 \subset I_2$ and $x=y$ on $I_1$; in case $I_1 \neq I_2$, we say that $y$ is a proper continuation of $x$. We say that $x:I \longrightarrow \R$ is a noncontinuable (or maximal) solution of (\ref{ivp}) if it has no proper continuation, and, in such a case, its domain $I$ is called maximal interval of definition.

Standard results guarantee that, in the conditions of Theorem \ref{pro1}, every solution can be continued as a solution of (\ref{ivp}) to a maximal interval of definition. Moreover, maximal intervals of definition are necessarily of the type $[t_0,T)$, for some $T>t_0$ and $T$ may be $+\infty$. Furthermore, see \cite[Theorem 1.3, Ch. 2]{codlev}, \cite[Theorem 4.1]{reid}, 
\cite[Pages 83--85]{sansone}, or \cite[Page 93]{wal}, the least and the greatest noncontinuable solutions exist.

For the sake of completeness we include our own proof of existence of the greatest noncontinuable solution. The existence of the least noncontinuable solution is identical and we omit it.

\begin{theorem}
In the conditions of Theorem \ref{pro1} there exist a noncontinuable greatest solution $x^*:[t_0,T) \longrightarrow \R$, where $T>t_0$ can be $+\infty$, in the sense that $x^*$ is a noncontinuable solution of (\ref{ivp}) and for any other solution $x:I \longrightarrow \R$ we have $x (t) \le x^*(t)$ for all $t \in I \cap [t_0,T)$.
\end{theorem}

\noindent
{\bf Proof.} Let ${\cal T}$ denote the set of times $\tau>t_0$ such that (\ref{ivp}) has the greatest solution on the interval $[t_0,\tau]$ in the following sense: there exists a solution $x_{\tau}:[t_0,\tau] \longrightarrow \R$ such that for any other solution $x:I \longrightarrow \R$ we have $x(t) \le x_{\tau}(t)$ for all $t \in I \cap [t_0,\tau]$.

To show that the set ${\cal T}$ is not empty, we go back to the proof of Theorem \ref{pro1} and we show that $t_0+L \in T$, where $L$ is defined after $M$ in condition (\ref{bded}). With the notation of that proof, we know that $x^*:[t_0,t_0+L] \longrightarrow \R$ is greater than any solution defined on $[t_0,t_0+L]$, so we have but to show that if $x:J \longrightarrow \R$ is another solution of (\ref{ivp}) with $J \subsetneq  [t_0,t_0+L]$ then $x \le x^*$ on $J$. First, we extend $x$ to a maximal interval of definition $[t_0,t_0+\hat L)$; for simplicity, we denote by $x$ both the former solution and its extension. If we prove that $\hat L>L$ then we have $x(t) \le x^*(t)$ for all $t \in [t_0,t_0+L]$ and, in particular, $x \le x^*$ on $J$. Assume, on the contrary, that $\hat L \le L$. This implies that $(t,x(t)) \in \mbox{Int}(R)$ for all $t \in [t_0,t_0+\hat L)$, and then $|x'(t)| \le M$
for all $t \in [t_0,t_0+\hat L)$ by virtue of (\ref{bded}).  Therefore $x$ would have a limit at $t_0+\hat L$ and it could be extended past $t_0+\hat L$, a contradiction\footnote{Notice that we have extended the information given in Theorem \ref{pro1}: the greatest solution in Theorem \ref{pro1} is greater than any other solution, no matter whether its domain is smaller than that of $x^*$. This fatct will be used again at the end of the present proof.}.  

\bigbreak

We define
$$T=\sup {\cal T} \in \R \cup \{+\infty\},$$
and we shall prove that the greatest noncontinuable solution of (\ref{ivp}) is the function $x^*:[t_0,T) \longrightarrow \R$ defined as follows: for each $t \in [t_0,T)$ we take any $\tau \in (t,T)\cap {\cal T}$ and we define $x^*(t)=x_{\tau}(t)$, where $x_{\tau}$ is the solution of (\ref{ivp}) which corresponds to $\tau$ by definition of ${\cal T}$.

\bigbreak

This definition does not depend on the choice of $\tau$, because for different values $\tau_1, \tau_2 \in {\cal  T}$ we have $x_{\tau_1}=x_{\tau_2}$ in the intersection of their domains. Therefore, $x^*$ is well defined on $[t_0,T)$ and $x^*$ is greater than any other other solution thanks to its definition.

Let us show that $x^*$ cannot be extended on the right of $T$. There is nothing to prove if $T=+\infty$, so we assume that $T<+\infty$. If we can find a proper continuation of $x^*$ then the point $(T,x^*(T^-))$ is an interior point of the domain of $f$ and then, by Theorem \ref{pro1}, the initial value problem
\begin{equation}
\label{ivp2}
x'=f(t,x), \quad t \ge T, \, \, x(T)=x^*(T^-),
\end{equation}
has the greatest solution $y:[T,T+\varepsilon] \longrightarrow \R$ for some $\varepsilon>0$. Therefore the solution $x=x^*$ on $[t_0,T)$ and $x=y$ on $[T,T+\varepsilon]$ is the greatest solution on $[t_0,T+\varepsilon]$, which implies that $T+\varepsilon \in {\cal T}$, a contradiction with $T=\sup {\cal T}$.

For the sake of completeness, let us prove that the solution $x=x(t)$ that we have just defined is, as we claimed, the greatest solution of (\ref{ivp}) on $[t_0,T+\varepsilon]$. If another solution $z:[t_0,\hat T] \longrightarrow \R$ satisfies $z(t_1)>x(t_1)$ for some $t_1 \in (t_0,\hat T)$, then we arrive at a contradiction. First, if $t_1 <T$ then we have that $z(t_1) >x^*(t_1)$, which contradicts with the definition of $x^*$. Second, if $t_1 \in [T,T+\varepsilon)$  and $z \le x^*$ on $[t_0,t_0+T)$, then the pointwise maximum of $z$ and $y$ on the interval $[T,T+\varepsilon]\cap [T, \hat T)$ defines a solution  of (\ref{ivp2}) which is somewhere strictly greater than $y$, a contradiction with the fact that $y$ is the greatest solution of (\ref{ivp2}).  \qed
 
\section{Differential inequalities}
The most important application of the existence of the greatest solution of (\ref{ivp}) reveals when we prove that it is greater than any lower solution. 

\begin{definition}
\label{lo}
A lower solution of (\ref{ivp}) is a function $\alpha:I \longrightarrow  \R$ defined on an interval $I$ which has $t_0$ as minimum, $\alpha(t_0) \le x_0$, $\alpha$ is differentiable on $I$ and
$$\alpha'(t) \le f(t,\alpha(t)) \quad \mbox{for all $t \in I$.}$$

An upper solution of (\ref{ivp}) is a function $\beta:I \longrightarrow  \R$ defined on an interval $I$ which has $t_0$ as minimum, $\beta(t_0) \ge x_0$, $\beta$ is differentiable on $I$ and
$$\beta'(t) \ge f(t,\beta(t)) \quad \mbox{for all $t \in I$.}$$
\end{definition}

In the sequel, by ``greatest solution" we mean ``greatest noncontinuable solution." Also, every result that we are going to prove for greatest solutions and lower solutions has its analogue with reversed inequalities for least solutions and upper solutions.

\bigbreak

The following nonlinear comparison result, which in essence we owe to Peano and Perron \cite{peano, perron}, can be easily proven by means of a today standard brilliant idea by Scorza Dragoni \cite{sc1, sc} concerning a modified problem. See also \cite{caf, go, sch}.

\begin{theorem}
\label{gron}
In the conditions of Theorem \ref{pro1}, let  $x^*:[t_0,T) \longrightarrow \R$ be the greatest solution of (\ref{ivp}).

If $\alpha:I \longrightarrow \R$ is a lower solution of (\ref{ivp}) then $\alpha(t)�\le x^*(t)$ for all $t \in I \cap [t_0,T)$.

In particular, the greatest solution being a lower solution itself, we have the formula
\begin{equation}
\label{gs}
x^*(t)=\max \{ \alpha(t) \, : \, \mbox{$\alpha$ lower solution of (\ref{ivp})} \}, \quad t \in [t_0,T).
\end{equation}

 \end{theorem}

\noindent
{\bf Proof.}  Reasoning by contradiction, assume that $\alpha(t_1)>x^*(t_1) $ for some $t_1 \in (t_0,T)$, then there exist $t_2 \in [t_0,t_1)$ and $t_3 \in (t_2,t_1]$ such that
\begin{equation}
\label{recont}
x^*(t_2)=\alpha(t_2) \quad \mbox{and} \quad x^*(t)<\alpha(t), \, \, t \in (t_2,t_3].
\end{equation}

Consider the initial value problem
$$x'=\left\{
\begin{array}{cl}
f(t,x(t)), & \mbox{if $x(t) \ge \alpha(t)$,}�\\
\\
f(t,\alpha(t)), & \mbox{if $x(t)<\alpha(t)$,}
\end{array}
\right\}, \, \, t \ge t_2, \, \, \, x(t_2)=\alpha(t_2).$$
By Peano's existence theorem, this problem has at least one solution $\varphi:[t_2,t_2+\varepsilon] \longrightarrow \R$, for some $\varepsilon>0$ such that $t_2+\varepsilon <t_1$. Let us prove that $\varphi \ge \alpha$ on $[t_2,t_2+\varepsilon]$: if not, we could find $\hat t_2$ and $\delta>0$ such that $\varphi(\hat t_2)=\alpha(\hat t_2)$ and
$$\varphi(t)< \alpha(t) \quad \mbox{for all $t \in (\hat t_2, \hat t_2+\delta],$}$$
but then for every $t \in (\hat t_2,\hat t_2+\delta]$ we have
$$\varphi(t)-\varphi(\hat t_2)=\int_{\hat t_2}^{t}f(s,\alpha(s)) \, ds \ge \alpha(t)-\alpha(\hat t_2),$$
a contradiction, so the proof that $\varphi \ge \alpha$ is complete.

Then the solution $x=x^*$ on $[t_0,t_2]$ and $x=\varphi$ on $[t_2,t_2+\varepsilon]$ is somewhere strictly greater than $x^*$, a contradiction. 
\qed

\begin{remark}
There is nothing to change in the previous proof if we use the following more general definition of lower solution, which allows corners in their graphs: an integral lower solution is a continuous function $\alpha:[t_0,T) \longrightarrow \R$ such that for all $t_1, t_2 \in [t_0,T)$ such that $t_1 \le t_2$ we have
$$\alpha(t_2)-\alpha(t_1) \le \int_{t_1}^{t_2}{f(s,\alpha(s)) \, ds}.$$
A well--known proof of existence of greatest solutions starts with (\ref{gs}) as a definition, see \cite{go}. That approach is rather more involved than ours in Theorem \ref{pro1}. 
\end{remark}

The classical differential form of Gronwall--Bellman's Lemma \cite{bellman, gronwall} is now immediate from Theorem \ref{gron}. We can even prove very easily the following more general result, which needs no assumption on the sign of the relevant functions.

\begin{corollary}
\label{corgron}
Let $a, \, b:[t_0,T)\longrightarrow \R$ be continuous functions and let $\alpha:[t_0,T) \longrightarrow \R$ be a differentiable function satisfying
$$\alpha'(t) \le a(t) \, \alpha(t) +b(t), \, \, t \in [t_0,T).$$
Then
$$\alpha(t) \le \alpha(t_0)\, {\rm exp} \left( \int_{t_0}^{t}{ a(s) \, ds} \right)+\int_{t_0}^t{ b(s) \, \mbox{\rm exp} \left( \int_{s}^{t}{a(r) \, dr}\right) ds},\, \, t \in [t_0,T).$$
\end{corollary}

\noindent
{\bf Proof.} Just notice that $\alpha$ is a lower solution of the linear problem
$$x'(t)=a(t) \, x(t)+b(t), \, \, \, t \ge t_0, \, \, \, \, x(t_0)=\alpha(t_0),$$
which has a unique\footnote{Remember that this can be proven without any uniqueness theorem: we just have to multiply the differential equation by ${\rm exp}(-\int_{t_0}^t{a(s) \, ds})$ and integrate.} solution given by
$$x^*(t)=\alpha(t_0)\, {\rm exp} \left( \int_{t_0}^{t}{ a(s) \, ds} \right)+\int_{t_0}^t{ b(s) \, \mbox{\rm exp} \left( \int_{s}^{t}{a(r) \, dr}\right) ds}.$$
\qed

%
%
%
%
Finally, we establish the integral form of Gronwall's inequality, see \cite[Chapter III, Theorem 1.1]{hart}. We need an extra assumption and, for completeness, we repeat the proof of Corollary 4.4, Chapter III, in \cite{hart}.
\begin{corollary}
\label{gronin}
Let $f:D \longrightarrow \R$ be continuous and nondecreasing in $x$ for each fixed $t$.

If $x^*:I=[t_0,T) \longrightarrow \R$ is the greatest solution of (\ref{ivp}) then for every continuous function $\alpha(t)$ the inequality
$$\alpha(t) \le x_0+\int_{t_0}^t{f(s, \alpha(s)) \, ds}, \, \, t\in I,$$
implies 
$$\alpha(t) \le x_0+\int_{t_0}^t{f(s,\alpha(s)) \, ds} \le x^*(t), \, \, t \in I.$$

In particular, we have the usual integral form of Gronwall--Bellman's inequality: if $\alpha, \, a, \, b:[t_0,T)\longrightarrow \R$ are continuous functions, where $T>t_0$ can be $+\infty$, $a(t) \ge 0$ for all $t \in [t_0,T)$, and
$$\alpha(t) \le x_0+\int_{t_0}^t{a(s) \, \alpha(s) \, ds}+\int_{t_0}^t{b(s) \, ds}, \, \, t\in [t_0,T),$$
then 
$$\alpha(t) \le x_0\, {\rm exp} \left( \int_{t_0}^{t}{ a(s) \, ds} \right)+\int_{t_0}^t{ b(s) \, \mbox{\rm exp} \left( \int_{s}^{t}{a(r) \, dr}\right) ds},\, \, t \in [t_0,T).$$

\end{corollary}

\noindent
{\bf Proof.} Let $\alpha(t)$ be as in the statement and define
$$\hat \alpha(t)=x_0+\int_{t_0}^t{f(s,\alpha(s)) \, ds}, \, \, t \in I.$$
We then have $\hat \alpha(t_0)=x_0$ and
$$\hat \alpha'(t)=f(t,\alpha(t))�\le f(t,\hat \alpha(t)), \, \, t \in I,$$
which implies $\hat \alpha \le x^*$ on $I$ by virtue of Theorem \ref{gron}.

\bigbreak

The proof of Gronwall--Bellman's inequality is trivial from the first part and we omit it. Note however that $a(t)$ must be nonnegative so that $f(t,x)=a(t)x+b(t)$ be nondecreasing in $x$.
\qed

\section{Peano's inequalities}
Peano's inequalities have to do with strict lower (upper) solutions and give finer information which, in particular, yields strict forms of Gronwall's inequality. Let us clarify that the denomination ``Peano's inequality" is not standard. We have adopted it for the sake of acknowledging Peano's ideas in \cite{peano}, when he proved that the least solution of (\ref{ivp}) is the pointwise supremum of all strict lower solutions. 
 
 \begin{theorem}
\label{peain}
(Peano's inequality for strict lower solutions) Assume that $f:D \longrightarrow \R$ is continuous and let  $x_*:[t_0,T) \longrightarrow \R$ be the least solution of (\ref{ivp}).

If $\alpha:[t_0,T) \longrightarrow \R$ is continuous and

\bigbreak
 
\hspace*{1cm} either $\alpha(t_0) <x_0$, $\alpha'$ exists in $(t_0,T)$, and
$$\alpha'(t) < f(t,\alpha(t)) \quad \mbox{for all $t \in (t_0,T)$,}$$

\medbreak

\hspace*{1cm} or $\alpha(t_0) = x_0$, $\alpha'$ exists in $[t_0,T)$, and
$$\alpha'(t) < f(t,\alpha(t)) \quad \mbox{for all $t \in [t_0,T)$,}$$
then $\alpha(t)<x_*(t)$ for all $t \in (t_0,T)$.
 \end{theorem}
 
 \noindent
 {\bf Proof.} Assume first that $\alpha(t_0)<x_0$; reasoning by contradiction, we assume that we can find some $t_1 \in (t_0,T)$ such that 
 \begin{equation}
 \label{econp}
\alpha(t_1)=x_*(t_1) \quad \mbox{and} \quad \alpha(t)<x_*(t) \quad \mbox{for all $t \in [t_0,t_1)$.}
 \end{equation}
 Then we have
 $$\alpha'(t_1)<f(t,\alpha(t_1))=f(t,x_*(t_1))=x_*'(t_1),$$
 a contradiction with (\ref{econp}).
 
 In case $\alpha(t_0)=x_0$ we have $\alpha'(t_0)<f(t,\alpha(t_0))=x_*'(t_0)$, hence $\alpha < x_*$ in $(t_0,t_0+\varepsilon)$ for some $\varepsilon>0$. In order to prove that $\alpha < x_*$ on the whole of $(t_0,T)$ it suffices to repeat the arguments in the first part with $t_0$ replaced by any $\hat t_0 \in (t_0,t_0+\varepsilon)$.\qed
 
 Notice that we cannot replace the assumptions on $\alpha$ by the following: ``$\alpha(t_0) \le x_0$ and $\alpha'(t)<f(t,\alpha(t))$ for $t�\in (t_0,T)$." Indeed, consider Peano's example
 $$x'=3x^{2/3},\,  \, \, t \ge 0, \, \,\, x(0)=0,$$
 which has $\alpha(t)=t^6$ as a strict lower solution on $(0,1/2)$, but not on $[0,1/2)$, and $\alpha$ is greater than the zero solution. Note, however, that $\alpha$ is less than the greatest solution, namely $x^*(t)=t^3$, which agrees with the result in Theorem \ref{gron}.

 \begin{remark}
 We could now go back and prove Theorem \ref{gron} by means of Theorem \ref{peain} with $x_*$ replaced by an arbitrary solution of (\ref{ivp}). This is another standard way: see, for instance, \cite[Theorem 11.6]{agor} and \cite[Lemma I.4]{halanay}.
 \end{remark}
 
The integral version of Theorem \ref{peain} can be obtained when $f(t,x)$ is strictly increasing in $x$ and adjusting the proof of Theorem \ref{gronin} accordingly. In doing so, we obtain a strict version of Gronwall--Bellman's inequality.
\begin{theorem}
\label{peanin}
Let $f:D \longrightarrow \R$ be continuous and strictly increasing in $x$ for each fixed $t$.

If $x_*:[t_0,T) \longrightarrow \R$ is the least solution of (\ref{ivp}) then for every continuous function $\alpha(t)$ the inequality
$$\alpha(t) < x_0+\int_{t_0}^t{f(s, \alpha(s)) \, ds}, \, \, t \in [t_0,T),$$
implies 
$$\alpha(t) < x_0+\int_{t_0}^t{f(s,\alpha(s)) \, ds} < x^*(t), \, \, t \in (t_0,T).$$

In particular, we have Gronwall's strict inequality: if $\alpha, \, a, \, b:[t_0,T)\longrightarrow \R$ are continuous functions, $a(t) > 0$ for all $t \in (t_0,T)$, and
$$\alpha(t) < x_0+\int_{t_0}^t{a(s) \, \alpha(s) \, ds}+\int_{t_0}^t{b(s) \, ds}, \, \, t\in [t_0,T),$$
then 
$$\alpha(t) < x_0\, {\rm exp} \left( \int_{t_0}^{t}{ a(s) \, ds} \right)+\int_{t_0}^t{ b(s) \, \mbox{\rm exp} \left( \int_{s}^{t}{a(r) \, dr}\right) ds},\, \, t \in [t_0,T).$$

\end{theorem}

\section{Maximal solution or greatest solution?}
It is customary in the literature to call maximal (minimal) solutions to what we mean in this paper by greatest (least) solutions. In our opinion, the adjective ``maximal" is inaccurate and misleading. It is surprising that at some moment the more precise denomination of ``maximum solution" that we find in \cite{codlev} was turned into ``maximal solution" and became extraordinarily popular! Our choice of calling greatest and least solutions is motivated from Carl and Heikkil\"a's book \cite{ch}.

\bigbreak

There are strong reasons for giving up saying maximal when we mean greatest:
\begin{enumerate}
\item The greatest solution is the maximum of the set of solutions equipped with the usual pointwise partial ordering, and not merely a maximal element of that set. Would you say 
``I have an animal" when you want to tell someone that you have a horse?
\item There exist differential problems which lack the greatest solution but do have maximal solutions in the right set theoretic sense; see \cite{h, fp}.
\item Sometimes, and specially in the Spanish mathematical literature, maximal solution means noncontinuable solution.
\end{enumerate}

\end{document}